\documentclass[11pt]{amsart}
\usepackage{latexsym,amssymb,amsmath}
\def\t{\tau}
\def\s{\sigma}
\def\overarrow{\vec}

\def\tbase{\text{Base}\,}

\def\BC{\mathbb C}

\def\BP{\mathbb P}
\def\pp#1{\mathbb P^{#1}}

\def\pp#1{{\mathbb P}^{#1}}
\def\tdim{\rm dim}
\def\hd{,...,}
\def\ww{\wedge}
\def\upperp{{}^\perp}

\def\cJ{{\mathcal J}}
\def\cE{{\mathcal E}}

\def\cS{{\mathfrak S}}

\def\cO{{\mathcal O}}
\def\CC{\mathbb C}

\def\11{\mathbf 1}
\def\PP{\mathbb P}

\def\FS{{\mathfrak S}}
\def\fsl{{\mathfrak {sl}}}

\def\fg{{\mathfrak g}}

\def\l{\lambda}
\def\a{\alpha}

\def\b{\beta}
\def\g{\gamma}
\def\s{\sigma}

\def\up#1{{}^{({#1})}}

\def\ot{{\mathord{\,\otimes }\,}}
\def\op{{\mathord{\,\oplus }\,}}
\def\otc{{\mathord{\otimes\cdots\otimes}\;}}

\def\lra{{\mathord{\;\longrightarrow\;}}}
\def\ra{{\mathord{\;\rightarrow\;}}}

\def\we{{\mathord{{\scriptstyle \wedge}}}}

\def\dim{{\rm dim}\;}
\def\La#1{\Lambda^{#1}}

\def\tdim{\text{dim}\,}

\def\tbase{\text{Base}\,}

\newcommand\rem{{\medskip\noindent {\em Remark}.}\hspace{2mm}}

\newtheorem{theo}{Theorem}

\newtheorem{theorem}{Theorem}[section]
\newtheorem{proposition}[theorem]{Proposition}
\newtheorem{lemma}[theorem]{Lemma}
\newtheorem{corollary}[theorem]{Corollary}

\theoremstyle{definition}
\newtheorem{definition}[theo]{Definition}

\newtheorem{example}[theo]{Example}

\theoremstyle{remark}

\begin{document}

\title{On the ideals of secant varieties of Segre varieties}
\author{J.M. Landsberg${}^1$
 \and L. Manivel}
\date{November 2003}
\footnote{Supported by NSF grant DMS-0305829}

\maketitle

\begin{abstract}
 We establish basic techniques for studying the ideals of
secant varieties of Segre varieties. We   solve a conjecture
of Garcia, Stillman and Sturmfels on the generators of the ideal of
the first secant variety in the case of three factors and solve the
conjecture set-theoretically for an arbitrary number of factors. We determine
the low degree
components of the ideals of secant varieties of small dimension in a few cases.
\end{abstract}

\section{Introduction}
Let $X^n\subset\BP V$ be a projective variety. Define $\s_r(X)$,
 the {\it variety of
secant $\pp{r-1}$'s to $X$} by
$$
\s_r(X)=\overline{ \cup_{x_1\hd x_r\in X}\BP_{x_1\hd x_r }}
$$
where $\BP_{x_1\hd x_r }\subset\BP V$ denotes the linear space
spanned by $x_1\hd x_r$ (usually a $\BP^{r-1}$).

\smallskip
Given $X\subset \BP V$ and $p\in\BP V$ define the {\it  essential $X$-rank
of $p$} (or {\it essential  rank of $p$} if $X$ is understood)  to be the smallest $r$
such that $p\in \s_r(X)$. 
(The essential rank is often called the {\it border rank} in
the computational complexity literature.) 
Similarly define the {\it rank} of $p$
to be the smallest $r$ such that there exist $r$ points on $X$,
$x_1\hd x_r$ such that $p\in \BP_{x_1\hd x_r}$. The essential
rank can be smaller than the rank, this phenomenon occurs
already for $X=v_3(\pp 1)$, the cubic curve, where the essential rank
of any point is at most two, but the rank of points on a tangent
line to $X$ (but not on $X$) is three.

\medskip

Given vector spaces $A_1\hd A_k$, one can form the {\it Segre product}
$X=Seg(\BP A_1\times \cdots\times \BP A_k)\subset \BP (A_1\ot\cdots\ot A_k)$.
When $k=2$, the Segre product is just the
projectivization of the space of rank one 
elements (matrices)
in $A_1\ot A_2$.  In this paper we study the ideals of the varieties
$\s_p(X)$.
The case $k=3$ is important in the study of computational
complexity as explained below. Many cases are important in the study
of Bayesian networks, as explained in \cite{GSS}. After presenting
some background information in \S2,  we establish
basic techniques for studying the problem for an arbitrary rational
homogeneous variety $X$ in \S 3.  In \S 4 we specialize to Segre products and
take advantage of Schur duality. We determine $I_3(\s_2(X))$ for any Segre product
 in Theorem \ref{cubicsprop}. We  prove that $I_3(\s_2(X))$ cuts out $\s_2(X)$
 set-theoretically in all cases and ideal theoretically when $k=3$, partially
 resolving Conjecture 21 of \cite{GSS}, see Theorem \ref{gssconj}.

\medskip

We present a deterministic algorithm to find 
  generators of the ideals of secant varieties
 of Segre varieties in \S 4. We carry this algorithm out in low degrees in \S 6.
  In particular,  we show there are no equations
in the ideal of $\s_6(\pp 3\times \pp 3\times \pp 3)$ in degree less than
nine. We plan to study higher degrees in a future paper.

One motivation for this paper is the following
question in computational complexity:
Let $A=(\BC^n\ot\BC^m)^*$, $B=(\BC^m\ot \BC^p)^*$, $C=\BC^n\ot\BC^p$.
The matrix multiplication operator
$M_{nmp}$ is an element of $A\ot B\ot C$.  In
standard coordinates $M_{nmp}$ is the sum of $nmp$ monomials. However,
already if one takes $m=n=p=2$, it is known that $\s_7(\pp 3\times\pp 3
\times \pp 3)=\pp {15}$ so the essential rank of $M_{222}$
is at most seven. Strassen showed \cite{Strassen} that in fact the
rank is at most seven by exhibiting an explicit expression of
$M_{222}$ as the sum of seven monomials, and moreover, he
proved that the essential rank is at least six by a specialization
argument.  The rank of $M_{222}$ was then shown to be seven in \cite{Winograd}.

To determine the essential rank of a point $p\in \BP V$ it is
sufficient to find   equation
 cutting out the varieties $\s_k(X)$ set theoretically
 and then to evaluate
the polynomials at $p$. Thus once one finds equations 
of $\s_6(Seg(\pp 3\times \pp 3\times \pp 3))$ one can determine
the essential rank of $M_{222}$.  

\medskip

\noindent{\it Acknowledgements} We thank Peter B\"urgisser for bringing the
border rank question to our attention and Bernd Sturmfels for comments on 
the exposition of a preliminary draft.

\section{Dimensions of secant varieties, especially homogenous ones}

Let $X\subset \BP V$ be a  projective variety.
An important fact about $\s_r(X)$ is Terracini's lemma, which implies
 that if $(x_1\hd x_r)$ is a general point  of $X\times X\times \cdots \times X$ then
 $$\hat T_{[\overarrow x_1+ ... +\overarrow x_r]}\s_r(X)
 =\hat T_{x_1}X+\cdots +\hat T_{x_r}X$$
 where $\overarrow x\in V$ denotes a point in the line $\hat x\subset V$
 corresponding to the point $x\in\BP V$ and $\hat T_pY\subset V$ denotes the affine tangent space to
 $Y$ at $p$, the cone over the embedded tangent projective
 space $\tilde T_pY\subset \BP V$.  
\medskip

If $X^n$ is smooth then the dimension of $\s_2(X)$ can be determined by
taking three derivatives at a general point   $x\in X$,
see \cite{GH}. In particular, if
the third fundamental form of $X$ at $x$,
$III_{X,x} $, is non-zero, then $\s_2(X)$ is of the expected dimension $2n+1$.

The third fundamental form calculation immediately 
implies that all homogeneously
embedded rational homogeneous varieties have  
$\s_2(X)$ of dimension $2\tdim X+1$ 
except for the following varieties
(embeddings are the  minimal homogeneous ones unless otherwise specified
and
if varieties occur more than one way we only list them once):  
$ 
G(2,n)=A_{n-1}/P_2$   the Grassmanian of
$2$-planes through  the origin in $\BC^n$, $Q^{2n-1}=D_n/P_1$,
$Q^{2n-2}=B_n/P_1$, the quadric hypersurfaces,
  $G_Q(2,2n)=D_n/P_2$, 
   $G_Q(2,2n+1)=B_n/P_2$, the Grassmanians of
$2$-planes throught the origin  
isotropic for a quadratic form, $v_{2}(\pp n)$ the quadratically embedded
Veronese,  $G_{\omega}(2,2n)=C_n/P_2$, the
the Grassmanians of
$2$-planes throught the origin  
isotropic for a symplectic form, $F_4/P_4$,
$G_2/P_1$, $E_6/P_1$, $E_6/P_2$, $E_7/P_1$,$E_7/P_7$, $E_8/P_8$.
$Seg(\pp k\times \pp l)=A_k/P_1\times A_l/P_1$.
Here we use the ordering of the roots as in \cite{bou}.

In all other cases the third fundamental form is easily seen to
be nonzero,  
see \cite{LM0}. In particular, for all triple and higher
Segre products, $\s_2(X)$ is nondegenerate.

In fact Lickteig and Strassen \cite{Lickteig,Strassenb} show many triple
Segre products have all secant varieties nondegenerate, in particular
for $Seg(\pp n\times\pp n\times \pp n)$ the filling secant variety
$\s_r$ is the expected number $r=\ulcorner n^3/(3n-2)\urcorner$ when $n>3$. In particular,
for $n=m^2$ we get roughly $m^4/3$ which is significantly greater than
$m^3$, so  in
higher dimensions matrix multiplication is far from being a generic
tensor. (Lickteig's proof is very simple and elegant - one first observes that
certain small cases, e.g., 
$\s_3(\pp 1\times \pp 1\times \pp n)$ fills   
and then one reduces to such cases by writing a larger vector space
as a sum of two dimensional spaces.)

\section{Ideals of secant varieties, especially homogeneous ones}

For $A\subset S^kV^*$ define $A\up p= (A\ot S^{p}V^*)\cap S^{p+k}V^*$,
the {\it $p$-th prolongation of $A$}. Let 
$$\tbase (A)=
\{ [v]\in \BP V\mid P(v)=0\ \forall P\in A\}.
$$
 Given a variety
$Z\subset\BP V$ we let $I(Z)\subset S^\bullet V^*$ denote
its ideal and $I_d(Z)=I(Z)\cap S^dV^*$.
We recall from \cite{LM0} that ideals of secant varieties satisfy the
{\it prolongation property}: 

\begin{lemma}[\cite{LM0}, Lemma 2.2]\label{prolonglemm} Let $A\subset S^2V^*$ be a system of quadrics with base locus $\tbase (A)\subset \PP V$. 
Then $$\tbase (A\up{k-1})\supseteq \s_k(\tbase (A)).$$
Moreover, if $\tbase (A)$ is linearly non-degenerate, then for $k\geq 2$, 
$I_k(\s_k(\tbase (A))=0$, and if $A=I_2(\tbase (A))$, then 
$I_{k+1}(\s_k(\tbase (A))=A\up k$.
\end{lemma}  

Geometrically,
$A\up{k-1}\supseteq \{ \frac{\partial P}{\partial v} \mid P\in A\up k\}
=I_{(\tbase A\up k)_{sing}}$ and $\s_{k-1}(X)\subseteq (\s_k(X))_{sing}$.

Usually $I(\s_k(X))$ is not generated in degree $k+1$. 
For example, consider the simplest intersection of quadrics, four points
in $\BP^2$. They generate six lines so $\s (X)$ is a hypersurface of
degree six. 

\begin{corollary}\label{emptycor} Let $X\subset \BP V$ be a variety with
$I(X)$ generated in degree $d$. Then for all $k\geq 0$,  $I_{d+k-2}(\s_{k}(X))=0$.
\end{corollary}

\smallskip

Given a variety $X\subset \BP V$,
a polynomial $P\in S^{d}V^*$, $d>k$, is in $I_d(\s_k(X))$ if and only if 
for any sequence of non negative integers $m_1,\ldots ,m_k$, with   $m_1+\cdots +m_k=d$, we have 
$P(v_1^{m_1},v_2^{m_2},\hdots  ,v_k^{m_k})=0$ for all $v_i\in \hat X$, where  
$v_i^{m_i}=v_i\circ\cdots\circ v_i$ ($m_i$ times).
Here we interpret $P(v_1^{m_1},v_2^{m_2},\hdots  ,v_k^{m_k})$ as the result of the successive contractions of $P$ 
by the tensors $v_i^{m_i} $.  

\smallskip

Now  consider the case where $X=G/P\subset \BP V_{\l}$ is a homogeneously
embedded rational homogeneous variety, i.e., 
the orbit of a highest weight
line.

By an unpublished theorem of Kostant,   $I_2(X)=(V_{2\l})\upperp \subset S^2V^*$ and $I(X)$
is generated in degree two. More generally,     $I_k(X)=(V_{k\l})\upperp
\subset S^kV^*$. We adopt the notation that if
$V=V_{\l}$, we write $V^k=V_{k\l}$.

Note that if   $W\subset S^dV^*$ is an
irreducible module, either all of $W$ or none of it is in
$I_d(\s_k(X))$.  These remarks imply:

\begin{proposition}\label{theidealisthis} 
 Let $X \subset \BP V$ be a rational homogenous   variety. Then a module
 $W\subset  S^dV^*$  is in
$I_d(\s_k(X))$ if and only if for all
integers $(a_1,\ldots ,a_p)$ such that $a_1+2a_2+\cdots +pa_p=d$ and $a_1+\cdots +a_p=k$,  the contraction map 
\begin{equation}\label{contractionmap}
W\ot S^{a_1}(V)\ot S^{a_2}(V^2)\otc S^{a_p}(V^p)\lra\CC
\end{equation}
is zero.
\end{proposition}

\begin{corollary}\label{turbostuffin}
Let $X=G/P\subset \BP V$ be a rational homogeneous variety. Then for all $d>0$, 
\begin{enumerate}
\item $I_d(\s_d(X))=0$, 
\item $I_{d+1}(\s_d(X))$ is the kernel of the contraction map $V^2\ot S^{d+1}V^*\ra S^{d-1}V^*$, 
\item let $W$ be an irreducible component of $S^dV^*$, and suppose that for all
integers $(a_1,\ldots ,a_p)$ such that $(a_1,\ldots ,a_p)$ such that $a_1+2a_2+\cdots +pa_p=d$ and $a_1+\cdots +a_p=k$,  
$W^*$ is not an irreducible component of $S^{a_1}(V)\ot S^{a_2}(V^2)\otc S^{a_p}(V^p)$. 
Then $W\subset I_d(\s_k(X))$. 
\end{enumerate}
\end{corollary}

\begin{proof}
(1) follows immediately from \ref{emptycor}, (2) from 
the remarks about the ideals of homogeneous varieties and
\ref{prolonglemm}. 
(3) follows from
  \eqref{contractionmap}
and Schur's lemma because if 
  an irreducible submodule
$W\subset S^dV^*$ does not belong to $I_d(\s_k(X))$, one of the contraction maps 
\eqref{contractionmap}
must be non-zero.
\end{proof}

\medskip\noindent {\bf Question}. {\it Is the ideal of the first secant variety $\s_2(X)$ of 
a rational homogeneous variety $X=G/P\subset\PP V$, generated by cubics 
(assuming its nonempty)? More generally,
when is the ideal of $\s_d(X)$ generated in degree $d+1$?}

\bigskip
If $X=G/P\subset \BP V$ is a Scorza variety, that is the set of rank one elements in
a simple Jordan algebra $\cJ$ (see \cite{Zak}), then $I_k(\s_{k-1}(X))$ is uniformly described as the
$k\times k$ minors in $\cJ$ and it generates $I(\s_{k-1}(X))$. More generally,
if $X=G/P$ is a sub-minuscule variety, that is, the set of tangent directions to lines
through a point of a compact Hermitian symmetric space (see \cite{LM0}), then
$I(\s_{k-1}(X))$ is generated in degree $k$ and there is a uniform description of 
$I_k$, see \cite{LM0}.

\smallskip
Let $\fg$ be a complex simple Lie algebra and consider 
the rational homogeneous variety $X_{ad}\subset\BP \fg$, the unique closed orbit of the 
corresponding adjoint Lie group. Then there are
universal modules $Y_k'\subset S^k\fg^*$, see \cite{LMuniv} and 
$Y_k'\subseteq I_k(\s_{k-1}(X_{ad}))$.

\smallskip
However, consider $X=Seg(\PP^2\times\PP^4\times\PP^6)
\subset \PP^{104}$. The expected dimension of $\s_8(X)$ is $8\times 12+7=103$,
thus if it is nondegenerate, $\s_8(X)$ is an invariant hypersurface. If its
degree is $d$, $S^d(A\ot B\ot C)$ must contain a one dimensional irreducible 
factor $(det A)^{\a}\ot (det B)^{\b}\ot (det C)^{\g}$. In particular $d$ is  divisible
by $3,5$ and $7$. Conclusion: either $\s_8(X)$ is degenerate, or it is a hypersurface 
of degree a multiple of $35$. This suggests that the degrees of the equations of the 
$\s_k(X)$ must  be much larger than $k$.  
 
\section{Schur duality and equations of Segre products}\label{schursect}

Let
$A_1\hd A_k$ be vector spaces and let $V=A_1\ot\cdots\ot A_k$. In order to determine the ideals
of secant varieties of Segre varieties $Seg(\BP A_1\times \cdots\times \BP A_k)\subset
\BP (A_1\ot \cdots \ot A_k)$, we need to understand
the decomposition of $S^dV^*$ into irreducible modules. We begin
by reviewing Schur duality (see,
e.g., \cite{FH}   for an introduction to Schur duality):

The irreducible representations
of the symmetric group $\cS_m$ are parametrized by the partitions of $m$. If $\pi$ is such
a partition, we let $[\pi]$ 
denote the corresponding $\cS_m$-module. For any vector 
space $V$, there is a natural action of $\cS_m$ on   $V^{\ot m}$ 
and we   define $S_{\pi}V$ the $\pi$-th {\it Schur power} of $V$ by
$$S_{\pi}V := {\rm Hom}_{\cS_m}([\pi],V^{\ot m}),$$ 
the $\cS_m$-equivariant linear maps from $[\pi]$ to $V^{\ot m}$. 
$S_{\pi}V$ is zero if $\pi$ has more
parts than the dimension of $V$, otherwise $S_{\pi}V$ is an irreducible $GL(V)$-module. 
Schur duality is the assertion that the tautological map 
$$\bigoplus_{|\pi |=m}[\pi]\ot S_{\pi}V \lra V^{\ot m}$$ 
is an isomorphism. 

For example, the trivial representation $[m]$ gives rise to
$S_mV=S^mV$ and the sign representation $[1\hd 1]$ gives rise
to $S_{1\hdots 1}V=\La m V$.

\begin{proposition} Let $A_1,\ldots ,A_k$ be vector spaces. Then 
$$S^m(A_1\otc A_k)= \bigoplus_{|\pi_1|=\cdots =|\pi_k|=m}([\pi_1]\otc [\pi_k])^{\cS_m}
S_{\pi_1}A_1\otc S_{\pi_k}A_k,$$
where $([\pi_1]\otc [\pi_k])^{\cS_m}$ denotes the space of $\cS_m$-invariants
(i.e., instances of the trivial representation) in the tensor
product. \end{proposition}

\proof  Apply Schur duality separately to each of $A_1,\ldots ,A_k$, take the 
tensor product of the corresponding isomorphisms, and compare with Schur duality 
for $A_1\otc A_k$. \qed 

\medskip 
Note that, since the representations of $\cS_m$ are self-dual, the dimension of  
$([\pi_1]\otc [\pi_k])^{\cS_m}$ is equal to the multiplicity of $[\pi_k]$ in the tensor 
product $[\pi_1]\otc [\pi_{k-1}]$. There is no general rule to compute such multiplicities,
but for small $m$ we can compute them using elementary character theory: if $\chi_{\pi}$ is
the character of $[\pi]$, then 
\begin{equation}\label{charcal}
\dim ([\pi_1]\otc [\pi_k])^{\cS_m} = \frac{1}{m!}\sum_{\s\in\cS_m}
\chi_{\pi_1}(\s)\cdots \chi_{\pi_k}(\s).
\end{equation}
  
\begin{proposition}\label{S3decomp} We have the following decomposition of
$S^3(A_1\ot \cdots\ot A_k)$ into irreducible $GL(A_1)\times \cdots
\times GL(A_k)$-modules:
\begin{multline*}
S^3(A_1\ot \cdots\ot A_k) =  \\
\bigoplus_{ \substack{|I|+|J|+|L|=k,\\ |J|>1}}
\frac{2^{j-1}-(-1)^{j-1}}3 S_3A_I\ot S_{21}A_J\ot S_{111}A_L 
\oplus
\bigoplus_{\substack{|I|+|L|=k,\\ |L|\,even}}S_3A_I\ot S_{111}A_L, 
\end{multline*}
where $I,J,L$ are multi-indices whose union is $1\hd k$, and we 
use the notation
$S_{\pi}A_I=\ot_{i\in I}S_{\pi}A_i$.

In particular, $S^3(A\ot B\ot C)=
S_3S_3S_3 \op S_3S_{21}S_{21}\op S_3S_{111}S_{111}\op
 S_{21}S_{21}S_{21}\op S_{21}S_{21}S_{111}$ and thus is multiplicity free.
Here $S_{\l}S_{\mu}S_{\nu}$ is to be read as
$S_{\l}A\ot S_{\mu}B\ot S_{\nu}C$  plus permutations giving rise to
distinct modules.
\end{proposition}

\begin{proof}
The irreducible representations of $\cS_3$ are the trivial representation $[3]$,  
the  sign representation $[111]$, and the natural two-dimensional representation $[21]$.
So we just need to compute the decomposition of $[21]^{\ot j}$ into irreducible components, 
which is a simple character computation. 

The symmetric group $\cS_3$ has three conjugacy classes of cardinality $1, 3, 2$, 
and the values of the irreducible characters on these classes are given by the following 
table:
$$\begin{array}{cccc}
{\rm class} & [Id] & [(12)] & [(123)] \\
\#    &  1  & 3 & 2 \\
\chi_3 & 1 & 1 & 1 \\
\chi_{21} & 2 & 0 & -1 \\
\chi_{111} & 1 & -1 & 1 
\end {array}$$
We calculate 
$$\langle \chi_{21}^j,\chi_3\rangle = \langle \chi_{21}^j,\chi_{111}\rangle = 
\frac{1}{6}(2^j+2(-1)^j)=\frac{1}{3}(2^{j-1}-(-1)^{j-1}),$$
where $\langle \chi,\chi'\rangle = \frac{1}{6}\sum_{\s\in\cS_3}\chi(\s)\chi'(\s)$ is the 
usual scalar product. The proposition follows. \end{proof}

The same type of computations lead to the following decomposition of the fourth 
symmetric power of a tensor product.

\begin{proposition}  We have the following decomposition of
$S^4(A\ot B\ot C)$ into irreducible $GL(A )\times GL(B)
\times GL(C)$-modules:
\begin{align}
S^4(A\ot B\ot C) = &\notag \; 
S_{4}S_{4}S_{4}\op S_{4}S_{31}S_{31}\op 
S_{4}S_{22}S_{22}\op S_{4}S_{211}S_{211}\op S_{4}S_{1111}S_{1111}\\
&\notag
\op S_{31}S_{31}S_{31}\op S_{31}S_{31}S_{22}\op S_{31}S_{31}S_{211}
\op S_{31}S_{22}S_{211}\op S_{31}S_{211}S_{211}\\
&\notag
\op S_{31}S_{211}S_{1111}\op S_{22}S_{22}S_{22}\op S_{22}S_{22}S_{1111}
\op S_{22}S_{211}S_{211}\op S_{211}S_{211}S_{211}.
\end{align}
Here $S_{\l}S_{\mu}S_{\nu}$ is to be read as
$S_{\l}A\ot S_{\mu}B\ot S_{\nu}C$  plus permutations giving rise to
distinct modules.
In particular, $S^4(A\ot B\ot C)$ is multiplicity free.
\end{proposition}

\rem In $S^5(A\ot B\ot C)$ all
submodules occuring have multiplicity one except for
$S_{311}S_{311}S_{221}$ which has multiplicity two.
For higher degrees there is a rapid growth in multiplicites.

\medskip
Now we try to determine   which modules are in the ideals of  the secant varieties
of $X=Seg(\PP A_1\times\cdots\times\PP A_k)$. We begin with some simple 
observations:

\begin{proposition}\label{elimprop}[Inheritance]
Suppose that an invariant $I$   of $[\l_1]\otc [\l_k]$ defines a nonzero embedding of 
$I$ into 
$S_{\l_1} A_1^*   \otc S_{\l_k}A_k^*\subset S^d(A_1\otc A_k)^*$.

  Then, for any vector spaces $A_1',\ldots,A_k'$ such that $\dim A_i'
\geq \dim A_i$ for all $i$, the image of the embedding of $(S_{\l_1}A'_1)^*\otc (S_{\l_k}A'_k)^*$ 
in $S^d(A'_1\otc A'_k)^*$ defined by $I$, is   in 
$I_d(\s_r(Seg(\BP A_1'\times \cdots \times \BP A_k')))$
if and only if the image of the embedding of 
$S_{\l_1} A_1^*   \otc S_{\l_k}A_k^*$ in $S^d(A_1\otc A_k)^*$
defined by $I$
  is 
in $I_d(\s_r(Seg(\BP A_1\times \cdots \times \BP A_k)))$.
\end{proposition}

\begin{proposition}\label{stuffin} Let $X^{(1)}=Seg(\PP A_2\times\cdots\times\PP A_k)$. Then
$$\begin{array}{rcl}
 I_d(\s_{d-1}(X))\cap (S^dA_1^*\ot S^d(A_2^*\otc A_k^*) ) & = & S^dA_1^*\ot 
I_d(\s_{d-1}(X^{(1)})).
\end{array}$$\end{proposition}

Say $S_{\pi_1}A_1\ot\cdots \ot S_{\pi_k}A_k\in S^d(A_1\ot\cdots \ot A_k)$.
An easy way to verify if it is in $I_d(\s_r(X))$ is if corollary
\ref{turbostuffin} (3) applies. This method is possible in low degrees but in higher degrees
multiplicities appear and the method becomes impossible to use.
Thus one either needs to understand the maps   (\ref{contractionmap}) or to write
down explicit polynomials and test them on $\s_r(X)$. One can either test a special polynomial
in a module on
a general point or test  a general polynomial in a module   at a special point.
The routines we used were more adapted to the first method.
We now describe two ways to explicitly write down polynomials.
The first has the advantage of producing the entire module,
the second of being quicker in producing a polynomial that is a highest weight vector.

Fix $\pi_1\hd\pi_k$ partitions of $d$. Compute $\tdim ([\pi_1]\ot \cdots \ot [\pi_k])^{\FS_d}$, call this
number $m$.

\bigskip\noindent{\bf ALGORITHM 1}

\begin{itemize}

\item Explicitly realize the representations $[\pi_j]$ of $\FS_d$.

\item Take independent elements $e_j\in [\pi_j]$ and average
$e_1\ot \cdots \ot e_k$ over $\FS_d$. The result is either a nontrivial
invariant $I$ or zero. Continue finding such elements $I$ until one has
$m$ independent such.

\item Choose   embeddings $S_{\pi_j}A_j\ra A_j^{\ot d}$, the images of
the invariants $I_s$, $1\leq s\leq m$ give the modules. 

\end{itemize}

\begin{example} 
Let $k=4$ and $d=3$. The space of invariants $([21]\ot [21]\ot [21]\ot [21])^{\cS_3}$
has dimension $2$. The representation $[21]$ of $\cS_3$ can be realized as the hyperplane $x_1+x_2+x_3=0$
 in $\CC_3$, and the action of $\cS_3$ is to permute  the coordinates 
$(x_1,x_2,x_3)$. A basis of $[21]$
 is given by $e=(1,-1,0)$ and $f=(0,1,-1)$. 
To obtain 
a basis of
the invariants of $[21]\ot [21]\ot [21]\ot [21]$ we   consider the natural 
basis 
of $[21]^{\ot 4}$
and apply the averaging operator over all translates by $\cS_3$. Applying this 
procedure to $eeee= e\ot e\ot e\ot e$ and $eeff= e\ot e\ot f\ot f$, we obtain the two invariants
$$\begin{array}{rcl}
I_1 & = & eeee+(e+f)(e+f)(e+f)(e+f)+ffff, \\
I_2 & =& 2eeee+eeef+eefe+efee+feee+3eeff+ \\
 & & \hspace{3cm}                                       +3ffee+fffe+ffef+feff+efff+2ffff.
\end{array}$$
Now   consider the space of $\cS_3$-equivariant  morphisms $u$ from 
$[21]$ to $V^{\ot 3}$, where $V$ is any vector space. Let $E=u(e)$.
Let $s_1$ denote the transposition $(12)$ and $s_2$ the transposition $(23)$.  Since $s_1(e)=-e$, 
  we get $s_1(E)=-E$. Since $f=s_2(e)-e$, 
we have $u(f)=s_2(E)-E$. And since $s_1(f)=e+f$, we must have $E-s_2(E)+s_1s_2(E)=0$.
The conclusion is that  $S_{21}V=Hom_{\cS_3}([21],V^{\ot 3})$ is isomorphic to the space 
of tensors $E\in V^{\ot 3}$ such that $s_1(E)=-E$ and $E-s_2(E)+s_1s_2(E)=0$.
\smallskip

Choose an invariant $J\in ([21]\ot [21]\ot [21]\ot [21])^{\cS_3}$ and consider
the embedding of 
$$S_{21}A_1\ot S_{21}A_2\ot S_{21}A_3\ot S_{21}A_4$$
 in $S^3(A_1\ot A_2
\ot A_3\ot A_4)$ that it defines. If $J=\a eeee+\cdots +\b ffff$ and $u_i\in S_{21}A_i$, 
the corresponding polynomial is defined by the equation 
$$\begin{array}{l}
P^J_{u_1,u_2,u_3,u_4}(a_1b_1c_1d_1, a_2b_2c_2d_2, a_3b_3c_3d_3)=\\
\hspace{4cm}
\a u_1(e)(a_1a_2a_3)u_2(e)(b_1b_2b_3)u_3(e)(c_1c_2c_3)u_4(e)(d_1d_2d_3)+\cdots +\\
\hspace{5cm} +\b u_1(f)(a_1a_2a_3)u_2(f)(b_1b_2b_3)u_3(f)(c_1c_2c_3)u_4(f)(d_1d_2d_3).
\end{array}$$
Now we evaluate $P^J_{u_1,u_2,u_3,u_4}$ on $\s_2(X)$, which means that we let 
$a_2=a_1$, $b_2=b_1$, $c_2=c_1$, $d_2=d_1$. Since $u_i(e)$ is skew-symmetric in its first
two arguments, its contribution will always be zero and we get
$$\begin{array}{l}
P^J_{u_1,u_2,u_3,u_4}(a_1b_1c_1d_1, a_2b_2c_2d_2, a_3b_3c_3d_3)= \\
\hspace{4cm}\b u_1(f)(a_1a_1a_3)u_2(f)(b_1b_1b_3)u_3(f)(c_1c_1c_3)u_4(f)(d_1d_1d_3),
\end{array}$$
so that the module defined by $J$ is in $I_3(\s_2(\BP A_1\times \cdots \times\BP A_4))$ if and only if $\b=0$. 
\end{example}

An immediate generalization of this argument leads to the following result:

\begin{proposition}
The space of modules in $I_3(Seg(\BP A_1\times \cdots\times \BP A_k))$
induced from  $[21]^{\ot k}$ is a codimension one subspace of
the modules in $S^3V^*$ induced from $[21]^{\ot k}$.
\end{proposition}

This proposition  allows
one  to determine the space of cubics vanishing on $\s_2(X)$. Indeed, every 
component of $S^3(A_1\otc A_k)^*$ involving a wedge power will do. Those involving 
a symmetric power are determined inductively by Proposition \ref{stuffin}. The only remaining 
term is $S_{21}A_1\otc S_{21}A_k$, whose multiplicity  equals $(2^{j-1}-(-1)^{j-1})/3$. 
The previous proposition means that the subspace vanishing on $\s_2(X)$ has multiplicity
 one less.  

\begin{theorem}\label{cubicsprop}
The space of cubics vanishing on the secant variety $\s_2(Seg(\BP A_1\times \cdots \times\BP A_k))$ is 
\begin{multline*}
 I_3(\s_2(Seg(\BP A_1\times \cdots \times\BP A_k)))  =  \bigoplus_{\substack{|I|+|J|+|L|=k,\\ |J|>1,\, |L|>0}}
 \frac{2^{j-1}-(-1)^{j-1}}{3} S_3A_I\ot S_{21}A_J\ot S_{111}A_L \\
\oplus \bigoplus_{\substack{ |I|+|J|=k,\\ |J|>1}}
(\frac{2^{j-1}-(-1)^{j-1}}3 -1)S_3A_I\ot S_{21}A_J
\oplus
\bigoplus_{\substack{|I|+|L|=k,\\ |L|>0\,even}} S_3A_I\ot S_{111}A_L.
\end{multline*}
\end{theorem}

\medskip
\begin{corollary}\label{cubetriple} Let $X= Seg(\BP A\times \BP B\times \BP C)$.
Then 
$$\begin{array}{rcl}
I_3(\s_2(X)) & = & (S^3A\ot \La 3 B\ot \La 3 C)^* \op
(\La 3A\ot S^3 B\ot \La 3 C)^* \op
 (\La 3A\ot \La 3 B\ot S^3 C)^* \\
 & & \op (S_{21}A\ot S_{21}B\ot \La 3 C)^*
\op (S_{21}A\ot \La 3 B\ot S_{21} C)^* \op 
(\La 3A\ot S_{21} B\ot S_{21} C)^*,
\end{array}$$
the space of $3\times 3$ minors of the three possible flattenings of $A\ot B\ot C$.
In particular, letting $a=\tdim A$, $b=\tdim B$, $c=\tdim C$, we have
$$\begin{aligned}
\tdim I_3(\s_2(X)) = &
  \frac{a b c}{72} \Big(-6 (a b+ac+bc) - 8 (a+ b+c) + 16 + 27 a b c  
- 5( a^2  b^2  c + a^2  b c^2 +  a b^2c^2)\\
\   & - 3 (a^2  b c + a b^2c+ a b c^2) + 5 a^2  b^2  c^2   
 + 2( a^2  b +a^2c +ab^2+ac^2+b^2c+bc^2)\\
 &
 +2( a^2 b^2+  2 a^2 c^2 +2 b^2 c^2)
  \Big).
\end{aligned}$$
\end{corollary}

In particular, we recover the data collected in \cite{GSS} and computed by {\sc Macaulay}
for the triple Segre products. 

\begin{corollary} Let  $A$, $B$ $C$, $D$ have dimensions $a, b, c, d$ respectively, we get 
the following number of cubic equations:
\begin{align}\notag
\tdim (I_3&(\s_2(\BP A\times \BP B\times \BP C\times \BP D))=\\
 \notag & \frac{abcd}{1296}\Big(368-72(a+b+c+d+ab+ac+ad+bc+bd+cd)-8(a^2+b^2+c^2+d^2)\\
&\notag -54(abc+abd+acd+bcd) +8( a^2b^2 + a^2c^2+ b^2c^2+ a^2d^2+ b^2d^2+ c^2d^2)+567abcd\\
 &\notag
+18(a^2bc+ab^2c+abc^2+a^2bd+ab^2d+a^2cd+ b^2cd+ ac^2d+ bc^2d+ abd^2+ acd^2+ bcd^2)\\
&\notag
-27abcd(a+b+c+d)+18(a^2b^2c+ a^2bc^2+ ab^2c^2+ b^2c^2d+a^2b^2d+ a^2c^2d+ b^2cd^2+ \\ 
&\notag
+ a^2bd^2+ ab^2d^2+ a^2cd^2+ac^2d^2+ bc^2d^2) 
+10(a^2b^2c^2 + b^2c^2d^2+ a^2b^2d^2+ a^2c^2d^2)\\
&\notag
-45abcd(cd+bd+ad+bc+ac+ab)-63abcd(abc+abd+adc+bcd)+143a^2b^2c^2d^2\Big).
\end{align}
\end{corollary}

This   recovers all the computations of cubic equations in \cite{GSS}.

\medskip

Before describing our second algorithm we do some preparation:

Fix a partition 
  $\pi =(p_1\hd p_f)$ of size $d=p_1\cdots +p_f$. For $\a_1\hd \a_f\in A^*$, let 
$$
F_{A }=(\a_1)^{\ot(p_1-p_2)}\ot (\a_1\ww \a_2)^{\ot(p_2-p_3)}\ot
\cdots 
(\a_1\ww\cdots\ww\a_{f-1})^{\ot p_{f-1}-p_f}\ot
(\a_1\ww\cdots\ww\a_{f})^{\ot p_f}\in (A^*)^{\ot d}.$$
When the $\a_1,\ldots ,\a_f$ vary, the subspace of $(A^*)^{\ot d}$ generated by 
the $F_A$'s is a copy of $S_{\pi}A^*$. In other words, we have defined an element 
of $Hom_{GL(A)}(S_{\pi}A^*,(A^*)^{\ot d})$, which is isomorphic to $ [\pi]$ by Schur duality. 

\bigskip\noindent{\bf ALGORITHM 2}

\begin{itemize}

\item For each $A_j$, choose a basis $\a^j_1\hd \a^j_{\tdim A_j}$, and it is better to choose a weight basis
with $\a^j_1$ a highest weight vector.
Continue the notation $m= \tdim ([\pi_1]\ot \cdots \ot [\pi_k])^{\FS_d}$.

\item Fix elements $\tau_1\hd \tau_k\in \FS_d$.
Let
\begin{equation}\label{Feqn}
\nonumber F(a_1^1\ot a^2_1\ot\cdots \ot a^k_1\hd a^1_d\ot a^2_d \ot\cdots \ot a^k_d)
=
F_{A_1}(a^1_{\tau_1(1)}\hd a^1_{\tau_1(d)})
\cdots
F_{A_k}(a^k_{\tau_k(1)}\hd a^k_{\tau_k(d)}).
\end{equation}

By construction $F\in (S_{\pi_1}A_1\ot\cdots \ot S_{\pi_k}A_k)^*$.

\item Now let
$$
P=
\sum_{\s\in \FS_d}
F_{A_1}(a^1_{\sigma\tau_1(1)}\hd a^k_{\sigma\tau_1(d)})
\cdots
F_{A_k}(a^k_{\sigma\tau_k(1)}\hd a^k_{\sigma\tau_k(d)}).
$$

By construction $P\in S^d(A_1 \ot\cdots\ot A_k)^*$,  and
$P\in (S_{\pi_1}A_1\ot\cdots \ot S_{\pi_k}A_k)^*$ as it is a sum of
terms in $(S_{\pi_1}A_1\ot\cdots \ot S_{\pi_k}A_k)^*$. 
Either $P$ is zero or it gives a nontrivial element of 
$(S_{\pi_1}A_1\ot\cdots \ot S_{\pi_k}A_k)^*\subset S^d(A_1 \ot\cdots\ot A_k)^*$.

\end{itemize}
Note that since we are choosing highest weight vectors, linear combinations
will also be highest weight vectors, thus we have a systematic way to look for
polynomials even when multiplicities occur.

\medskip

 In practice we implemented the algorithm  in two parts as follows:

\medskip\noindent{\bf Input}:
\begin{itemize}
\item  $k$, the number of vector spaces; 
\item $d_1\hd d_k$; their
dimensions; 

\item $d$, the degree of the polynomial to be constructed; 

\item $\pi_1\hd \pi_k$,
partitions of $d$
 
 \end{itemize}
 
{\sc Part one}: Finding the polynomials.
 
 \begin{enumerate}
 
 \item Calculate $m=dim ([\pi_1] \ot \cdots \ot [\pi_k])^{\FS_d}$ via a character
 calculation as in \eqref{charcal}. 
 
 \item  Choose a collection of permutations $T_1=(\t_1\hd \t_k)$ with $\t_j\in \FS_d$
 (without loss of generality, take $\t_1=Id$). Write out $F^{T_1}$ as in
 \eqref{Feqn} above and then average over $\FS_d$ to obtain a polynomial $P^{T_1}$
 as above.
 
 \item Test if $P^{T_1}$ is identically zero either by a symbolic calculation or
 by evaluating it at a randomly chosen point. If it is zero return to step 2.
 
 \item Repeat steps 2 and 3 for collections of permutations $T_2\hd T_m$,  only
 when repeating step 3, not only   test if the polynomial is nonzero, but
 also test that it is linearly independent from the polynomials already constructed.
 
 \end{enumerate}
 
 \medskip\noindent{\bf Output}:   a basis of highest weight vectors for  the isotypic
submodule of copies of  $(S_{\pi_1}A_1\ot\cdots\ot S_{\pi_k}A_k)^*$ inside $S^d(A_1\ot\cdots\ot A_k)^*$.
 
 \bigskip
 
{\sc Part two}: Testing if any modules are in the ideal.
 
\medskip\noindent{\bf Input}:
\begin{itemize}
\item   The polynomials $P^{T_1}\hd P^{T_m}$ constructed in part one.

\item  $p$: where we will test for generators of $I_d(\s_p(\BP A_1\times\cdots\times \BP A_k))$.

 \end{itemize}

 \begin{enumerate}
 
 \item Write $d= up+r$, with $u,r$ nonegative integers and $r<p$.
  Let $P= c_1P^{T_1}+\cdots +c_mP^{T_m}$ where the $c_j$'s are variables.
Pick $p$ vectors in each $A_i$ at random,
 $a^i_1\hd a^i_p$.
Considering $P$ as a multi-linear form,
let $\overline a_j= a^1_j\ot\cdots\ot a^k_j$
evaluate 
$$P( \overline a_1\hd \overline a_p, 
\overline a_1\hd \overline a_p \hd \overline a_1\hd \overline a_p \hd 
\overline a_1\hd \overline a_r).$$

\item 
Now pick $m-1$ more such sets of vectors and solve for
the $c_j$'s.

\item 
For simplicity, say there is a unique solution, test on
one more set of vectors using $P$ with the $c_j$'s replaced by
their solution values.
If one gets zero, one has a good candidate.

Warning: this is just one of many tests to perform to see if a
candidate is in the ideal - we begin with this one only because
in practice it has been quite useful. Hence the next step:

\item Now   test $P$   on all possible ways of choosing the last $r$
vectors from the set of first $p$ vectors (e.g., one needs to test the
possibility of the first vector occuring $r$ times instead of $r$ different
vectors etc...). Ideally do this symbolically, but one gets an answer with very
high probability by testing at random points.

 \end{enumerate}

\noindent{\bf Output}: Either ruling out
 the modules $(S_{\pi_1}A_1\ot\cdots\ot  S_{\pi_k}A_k)^*
\subset S^d(A_1\ot\cdots\ot A_k)^*$ from being in $I_d(\s_p(X)$ or
determination of an explicit copy of 
$(S_{\pi_1}A_1\ot\cdots\ot  S_{\pi_k}A_k)^*
\subset S^d(A_1\ot\cdots\ot A_k)^*$ that is  in $I_d(\s_p(X))$
described by its   highest weight vector.

\medskip

Here are some   examples:

\begin{example}\label{211211211ex}
Consider $S_{211}A\ot S_{211}B\ot S_{211}C\subset S^4(A\ot B\ot C)$. Without loss of generality
assume $\tdim A=\tdim B=\tdim C=3$.
Take
$$
F(a_1b_1c_1\hd a_4b_4c_4)=
\a^1(a_1)det(a_2,a_3,a_4)
\b^1(b_2)det(b_1,b_3,b_4)
\g^1(c_3)det(c_1,c_2,c_4)
$$
and let $P$ be the corresponding polynomial. A simple
evaulation at a random point shows $P$ is not identically zero.
(Compare with taking all permutations $\tau$ to be the identity, then the average
over $\FS_4$ is indeed zero.)
\end{example}

\begin{example}\label{333333333ex}
Consider $S_{333}A\ot S_{333}B\ot S_{333}C$. Without loss of generality take
$\tdim A=\tdim B=\tdim C=3$.
We take
\begin{align*}
F=&det(a_1,a_2,a_3)det(a_4,a_5,a_6)det(a_7,a_8,a_9)
det(b_2,b_3,b_4)det(b_5,b_6,b_7)\\
&
det(b_1,b_8,b_9)det(c_3,c_4,c_5)det(c_6,c_7,c_8)det(c_1,c_2,c_9)
\end{align*}
Here it is more delicate to see the corresponding polynomial $P$ is not identically
zero because there will be terms that appear several times. One needs to check that
they do not have signs cancelling. (For example, had we had any pair of indices occuring
three times in a determinant, the corresponding polynomial would be zero because
the transposition of the indices would produce the same terms with opposite signs.)
One can also verify with Maple that the corresponding polynomial is nonzero.
\end{example}

\begin{example} \label{3213213111}
Consider $S_{321}A\ot S_{321}B\ot S_{3111}C\subset S^6(A\ot B\ot C)$ which occurs with multiplicity four.
Let
\begin{align*}
F_{\tau,\mu}=& \a_1\ww\a_2\ww\a_3(a_1,a_2,a_3)*
 \a_1\ww\a_2(a_5,a_6) *\a_1(a_4)\\
&
\b_1\ww\b_2\ww\b_3(b_{\t (1)} ,b_{\t (2)} ,b_{\t (3)} )
 *\b_1\ww\b_2(b_{\t (4)} ,b_{\t (5)})*\b_1(b_{\t (6)} ) \\
&
\g_1\ww\g_2\ww\g_3\ww\g_4 (c_{\mu (1)} ,c_{\mu (2)} ,c_{\mu (3)},c_{\mu (4)} )
*\g_1 (c_{\mu(5)} )*\g_1(c_{\mu (6)} )
\end{align*}
Now we take the following permutations:
\begin{align*} \t_1 &=\begin{pmatrix} 
1 &2&3&4&5&6\\ 
3&4&5&1&2&6\end{pmatrix}\ \qquad \mu_1
=\begin{pmatrix} 1 &2&3&4&5&6\\ 1&4&5&6&2&3\end{pmatrix}\\
\t_2 &
=\begin{pmatrix} 1 &2&3&4&5&6\\ 3&4&5&1&2&6\end{pmatrix}\ \qquad \mu_2
=\begin{pmatrix} 1 &2&3&4&5&6\\ 2&3&5&6&1&4\end{pmatrix}\\ 
\t_3 &
=\begin{pmatrix} 1 &2&3&4&5&6\\ 3&4&5&1&2&6\end{pmatrix}\ \qquad \mu_3
=\begin{pmatrix} 1 &2&3&4&5&6\\ 2&3&4&5&1&6\end{pmatrix}\\  
\t_4 &
=\begin{pmatrix} 1 &2&3&4&5&6\\ 3&4&6&1&2&5\end{pmatrix}\ \qquad \mu_4
=\begin{pmatrix} 1 &2&3&4&5&6\\ 2&3&4&5&1&6\end{pmatrix}
\end{align*} 

The resulting four polynomials, call them $P_1\hd P_4$ are linearly
independent. We verified  this by evaluating them first at four random points to
determine a unique possible linear combination that is zero, and then evaluated
this linear combination at a fifth random point - one does not obtain zero.
\end{example}

\rem When the $A_i$'s have the same dimension $k$,
$S^k(A_1\otc A_k)^*$ contains a copy of $\fsl (A_1)\otc \fsl (A_k)$, 
with an embedding given by the formula
$$\begin{array}{rcl}
P_{X_1,\ldots X_k}(a_1^1\otc a_k^1,\ldots,a_1^k\otc a_k^k)
 & = & (X_1a_1^1\wedge a_1^2\wedge\cdots\wedge a_1^k)\cdots 
(X_ka_k^1\wedge a_k^2\wedge\cdots\wedge a_k^k) \\
 & &  \quad + {\rm symmetric\; terms}.
\end{array}$$
All such polynomials vanish on $\s_{k-3}(X)$, but not on $\s_{k-2}(X)$.

\section{Flattenings and the GSS conjecture}

Let  $X=Seg(\PP A_1\times\cdots\times\PP A_k)\subset\PP (A_1\otc A_k)$. 
A family of degree $d+1$ equations for $\s_d(X)$ is given by the 
{\it flattenings}, discussed in \cite{GSS}.    

\begin{definition}
Given $V=A_1\ot\cdots \ot A_k$, a {\it flattening} of $V$ is a decomposition
$$
V=(A_{i_1}\ot\cdots \ot A_{i_q})\ot (A_{j_1}\ot\cdots \ot A_{j_{k-q}})= A_I\ot A_J
$$
where $I+J=\{1\hd k\}$ is a partition of $\{ 1\hd k\}$ into two subsets.
\end{definition}
 
Since $X\subset Seg(\BP A_I\times \BP A_J)$, $\s_k(X)\subseteq \s_k(Seg(\BP A_I\times \BP A_J))$
and thus the $(d+1)\times (d+1)$  minors 
of flattenings always vanish on $\s_d(X)$, i.e.
$$\wedge^{d+1}(A_{i_1}\ot\cdots \ot A_{i_q})^*
\ot \wedge^{d+1}(A_{j_1}\ot\cdots \ot A_{j_{k-q}})^*\subset I_{d+1}(\s_d(X)).$$
In \cite{GSS} it was conjectured that $I (\s_2(X))$ is generated
by the $3\times 3$ minors of flattenings, i.e., that $\s_2(X)$ is  
intersection
 as a scheme  of the varieties $\s_2(\BP A_I\times \BP A_J)$. We will prove this for $k=3$ below. 
For $k>3$ we have the following  partial result which implies that that $\s_2(X)$ is  
intersection
 as a set of the varieties $\s_2(\BP A_I\times \BP A_J)$. 

\begin{theorem}\label{gssconj}Let $X=Seg(\PP A_1\times\cdots\times
\PP A_k)\subset\PP (A_1\otc A_k)$ be a Segre product of projective spaces. 
 \begin{itemize}
\item The first secant variety $\s_2(X)$  is defined set theoretically by the 
$3\times 3$ minors of flattenings.

\item $I_3(\s_2(X))$ is spanned  by the 
$3\times 3$ minors of flattenings.
\end{itemize}
\end{theorem}

The corresponding modules were described explicitly in  Theorem \ref{cubicsprop}.

\proof Let $T\in A_1\otc A_k$ be a tensor on which the $3\times 3$ minors of 
flattenings all vanish. This means that $T$ has rank two at most, when considered 
has a tensor of $A_I\ot A_J$, where $A_I=A_{i_1}\ot\cdots \ot A_{i_q}$ and 
$A_J=A_{j_1}\ot\cdots \ot A_{j_{k-q}}$, with $V=A_I\ot A_J$ any flattening. 
Applying this to the case where $\# I=1$, we see that we can find two dimensional
subsets $A_i'\subset A_i$ such that $T\in A_1'\otc A_k'$. In other words, we may
and will suppose that $\dim A_i=2$ for all $i$. 

\smallskip
Now take $I=\{1,2\}$. We can decompose our tensor as $T=M\ot S+M'\ot S'$, 
where $M, M'\in A_1\ot A_2$. We can identify $A_1$ with the dual of $A_2$ and 
consider $M$ and $M'$ as endomorphisms of $A_2$. Suppose that one of them has rank
two. We can adapt our basis so that $M$, for example, is the identity and $M'$ is in 
Jordan canonical form. Generically, $M'$ will be diagonalizable and we can rewrite
our tensor as 
$$T=a_1\ot a_2\ot C+a_1'\ot a_2'\ot C'.$$
If $a_1$ and $a_1'$, or $a_2$ and $a_2'$, are proportional, $T$ can be factored as 
$a_1\ot U$ and we are reduced to the case of $k-1$ factors. So we can suppose that 
$(a_1,a_1')$ is a basis of $A_1$, and $(a_2,a_2')$ a basis of $A_2$. Then we consider
$C$ and $C'$ as map from $(A_4\otc A_k)^*$ to $A_3$ and apply our hypothesis to the
set of indices $I=\{1,3\}$. The conclusion is that $a_1\ot C(t)$ and $a_1'\ot C'(t)$ 
belong to a fixed two-dimensional subset of $A_1\ot A_3$, as $t$ varies in  $(A_4\otc A_k)^*$. 
Since $a_1$ and $a_1'$ are independant, this implies that $C$ and $C'$ have rank one. 
But the same conclusion holds if we replace $I=\{1,3\}$ by any $I=\{1,j\}$, $j\ge 3$, 
and this means that we can decompose $C=a_3\otc a_k$ and $C'=a_3'\otc a_k'$. Thus 
$T$ belongs to the secant variety $\s_2(X)$.

Suppose now that $M'$ is not diagonalizable. Then we can find bases $(a_1,a_1')$ 
of $A_1$, and $(a_2,a_2')$ of $A_2$, such that we can decompose $T$ as
$$T=(a_1'\ot a_2+a_1\ot a_2')\ot C+a_1\ot a_2\ot C'.$$
We shall prove by induction on $j\ge 2$ that we can decompose $T$ further as 
$$T=(a_1'\ot a_2 \otc a_j+\cdots +a_1\otc a_{j-1}\ot  a_j')\ot C_j+a_1\otc a_j\ot C_j',$$
for some $C_j,C_j'\in A_{j+1}\otc A_k$. 
As in the previous case, we consider $C_j$ and $C_j'$ as morphisms from $(A_{j+2}\otc 
A_k)^*$ to $A_{j+1}$. Then $a_2 \otc a_j\ot C_j(t)$ and $(a_2'\otc a_j+\cdots+
a_2\otc a_j')\ot C_j(t)+a_2\otc a_j\ot C_j'(t)$ belong to a fixed two dimensional space
$V_j$ as $t$ varies. This implies that $C_j$ has rank one, we write it as $C_j=a_{j+1}\ot C_{j+1}$. 
Then $V_j$ contains the tensors $a_2\otc a_{j+1}$, $(a_2'\otc a_j+\cdots +a_2\otc a_j')\ot a_{j+1}
+a_2\otc a_j\ot C_j'(t_0)$ for $C_{j+1}(t_0)=1$, and $a_2\otc a_j\ot C_j'(t)$ for $t$ in the 
kernel of $C$. But the first two vectors are already independant, so that those of the third 
type must be proportional to the first one, 
which means that $C_j'$ maps the kernel of $C_j$ its image. But this
means that we can decompose $C_j'=a_{j+1}\ot C_{j+1}'+a_{j+1}'\ot C_{j+1}$ for some 
$a_{j+1}'\in A_{j+1}$ and $C_{j+1}\in A_{j+2}\otc A_k$. This concludes the induction. 

\smallskip When $j=k-1$, we finally get a decomposition of $T$ as 
$$T=a_1'\ot a_2 \otc a_k+\cdots +a_1\otc a_{k-1}\ot  a_k'+a_1\otc a_k.$$
We conclude that $T$ belongs to the (affine) tangent space of $X$ at the point 
$a_1\otc a_k.$ In particular, $T$ belongs to the tangential variety of $X$, which is contained 
in the secant variety $\s_2(X)$.\qed \medskip

\begin{theorem} Let $X=Seg(\PP A\times\PP B\times \PP C)\subset \PP(A\ot B\ot C)$ be a triple Segre 
product. Then the ideal of the secant variety  $\s_2(X)$ is generated by cubics.
\end{theorem}

\proof Let $\hat{\s}\subset A\ot B\ot C$ denote the cone over  
$\s_2(\BP A\times\BP B\times \BP C)$.
Let $G_2(A)$ denote the Grassmanian of two-planes in $A$.
 Consider the quasiprojective variety
$Y = G_2(A)\times G_2(B)\times G_2(C)\times A\ot B\ot C$, and denote by $p$
and $\pi$ its projections to $G_2(A)\times G_2(B)\times G_2(C)$ and $A\ot B\ot C$. 
Let  $T_A$ denote the tautological rank two vector bundle on $G_2(A)$, 
and let $ E$ denote
the vector bundle $A\ot B\ot C/T_A\ot T_B\ot T_C$ on $G_2(A)\times G_2(B)\times G_2(C)$.
The pull-back $p^* E$ has a canonical section $s$, defined by 
$$s (U_A,U_B,U_C,t)=t \quad {\rm mod}\; U_A\ot U_B\ot U_C.$$
Let $\tilde{\s}$ denote the zero-locus of this section. 

\begin{lemma} The zero-locus $\tilde{\s}$ is a vector bundle over 
$G_2(A)\times G_2(B)\times G_2(C)$, in particular it is a smooth variety. 
Its image under $\pi$ is $\hat{\s}$, and the restriction map $\pi_{|\tilde{\s}} :
\tilde{\s}\ra\hat{\s}$ is a resolution of singularities. 
\end{lemma}

\proof The first assertion is clear. The second one is an immediate consequence of 
the fact that the secant variety of $\PP^1\times\PP^1\times\PP^1\subset\PP^7$
is non degenerate, i.e.,  equal to $\PP^7$. \qed 

\medskip
Consider the Koszul complex of the section $s$: this is a minimal free resolution of 
the structure sheaf of $\tilde{\s}$. We want to push it down to $A\ot B\ot C$ to get some 
information on the minimal resolution of $\hat{\s}$. For this we use the spectral 
sequence 
$$\cE^{p,q}_1 = R^q\pi_*p^*(\wedge^{-p} E^*) \Longrightarrow R^{p+q}\pi_*\cO_{\tilde{\s}}.$$

\begin{lemma} We have $R^q\pi_*\cO_{\tilde{\s}}=0$ for $q>0$ and $\pi_*\cO_{\tilde{\s}}=
\cO_{\hat{\s}}$. In particular, $\hat{\s}$ has rational singularities.
\end{lemma}

\proof The fibers of $\pi$ are isomorphic to $G_2(A)\times G_2(B)\times G_2(C)$, 
and the vector bundle $p^* E$ is a pull-back from that product. This reduces the
problem to the computation of the cohomology of $E^*$ and its exterior powers.
Let $V$ denote the trivial bundle whose fiber is isomorphic to $A\ot B\ot C$
and let $T=T_A\ot T_B\ot T_C$.
For each integer $r$, we have an exact sequence
$$0\ra\wedge^rE^*\ra \wedge^rV^* \ra \wedge^{r-1}V^*\ot T^* \ra 
\wedge^{r-2}V^*\ot S^2T^* \ra \cdots\ra S^rT^* \ra 0.$$
By Bott's theorem, the vector bundles $S^kT^*$ are acyclic. Since the previous
resolution of $\wedge^rE^*$ is of length $r+1$, this implies that 
$$\begin{array}{rcl}
H^q(G_2(A)\times G_2(B)\times G_2(C),\wedge^rE^*) & = & 0 \quad {\rm for}\; q>r, \\
H^r(G_2(A)\times G_2(B)\times G_2(C),\wedge^rE^*) & = & 
coker \Big(V^*\ot H^0(S^{r-1}T^*)\lra H^0(S^rT^*)\Big).
\end{array}$$
The first claim implies that in the spectral sequence ($r=-p$ !), $\cE^{p,q}_1=0$
for $p+q>0$, hence $R^k\pi_*\cO_{\tilde{\s}}=0$ for $k>0$. 

To prove that $\pi_*\cO_{\tilde{\s}}=\cO_{\hat{\s}}$, we need to check that 
$H^r(G_2(A)\times G_2(B)\times G_2(C),\wedge^rE^*)=0$ for $r>0$. But note that 
if 
\begin{eqnarray}
\nonumber  S^rV^* & = & \oplus_{\l,\mu,\nu}c_{\l,\mu,\nu}S_{\l}A^*\ot S_{\mu}B^*\ot 
S_{\nu}C^*, \\
\nonumber  {\rm then} \; 
S^rT^* & = & \oplus_{l(\l),l(\mu),l(\nu)\le 2}
c_{\l,\mu,\nu}S_{\l}T_A^*\ot S_{\mu}T_B^*\ot S_{\nu}T_C^*, \\
\nonumber  {\rm thus} \;
H^0(S^rT^*) & = & \oplus_{l(\l),l(\mu),l(\nu)\le 2}
c_{\l,\mu,\nu}S_{\l}A^*\ot S_{\mu}B^*\ot S_{\nu}C^*.
\end{eqnarray}
Therefore, if $r>0$, the surjectivity of the map $V^*\ot H^0(S^{r-1}T^*)\lra H^0(S^rT^*)$
is an immediate consequence of the surjectivity of $V^*\ot S^{r-1}V^*\lra S^rV^*$. \qed

\medskip We are now in position to apply Theorem (5.1.3) of \cite{weyman}, 
following which  the vector bundles $R^q\pi_*p^*(\wedge^{-p} E^*)$ can be organized 
into a resolution of $\cO_{\hat{\s}}$. In particular, the cohomology groups
$H^{r-1}(\wedge^rE^*)$ appear as degree $r$ equations of $\hat{\s}$. Thus, 
if we can prove that these groups vanish for $r\neq 3$, we'll get that the ideal
of  $\hat{\s}$ is generated by cubics. 

Using the previous resolution of $\wedge^rE^*$, 
we see that $H^{r-1}(\wedge^rE^*)$ is the homology group of the complex on the
first line of the diagram 
$$\begin{array}{ccccc}
\wedge^2V^*\ot H^0(S^{r-2}T^*) & \lra & V^*\ot H^0(S^{r-1}T^*) & \lra & H^0(S^rT^*) \\
\uparrow & & \uparrow & & \uparrow \\
\wedge^2V^*\ot S^{r-2}V^* & \lra & V^*\ot S^{r-1}V^* & \lra & S^rV^*
\end{array}$$
The complex on the lowest line is a Koszul complex. It is exact, and surjects onto
the complex we are interested in. For $r=1$ or $r=2$ we get the same complexes, hence
$H^0(E^*)=H^1(\wedge^2E^*)=0$. For $r=3$ only the rightmost terms are different,  
and $H^0(S^3T^*)$ is the sum of components in $S^3V^*$ without terms of length three. 
The other components  give $H^2(\wedge^3E^*)$.

Next we must prove that $H^{r-1}(\wedge^rE^*)=0$ for $r\ge 4$. First observe that 
the components of $V^*\ot H^0(S^{r-1}T^*)$ have length three at most. Those of length
at most two on each factor map to $H^0(S^{r}T^*)$ as they do in the Koszul complex,
which is exact: this takes care of that kind of terms. Now consider an isotypical
component $D$ inside $V^*\ot H^0(S^{r-1}T^*)$ with length three, say, on $A$. 
This component maps to zero in  $H^0(S^{r}T^*)$,
and we must check that it belongs to the image of $\wedge^2V^*\ot H^0(S^{r-2}T^*)$. 
But we know that the secant variety $Seg(\PP A\times\PP (B\ot C))$ is cut out by 
cubics, and this implies that the corresponding complex
$$\wedge^2V^*\ot H^0(S^{r-2}U^*)  \lra  V^*\ot H^0(S^{r-1}U^*)  \lra  H^0(S^rU^*)$$
is exact. Here $U$ is the tautological vector bundle $T_A\ot T_{B\ot C}$ on $G_2(A)
\times G_2(B\ot C)$, so that 
$$H^0(S^rU^*) = \oplus_{l(\l)\le 2, \mu, \nu}
c_{\l,\mu,\nu}S_{\l}A^*\ot S_{\mu}B^*\ot S_{\nu}C^*.$$
Therefore, we  see our component $D$ inside $V^*\ot H^0(S^{r-1}U^*)$, and for the same
reason as before, it maps to zero in $H^0(S^{r}U^*)$. So it must belong to 
the image of $\wedge^2V^*\ot H^0(S^{r-2}U^*)$. We must check that in fact, it only comes 
from components of $H^0(S^{r-2}U^*)$ with length at most two on each factor. But this 
is clear, because the contraction map factors as 
$$\wedge^2V^*\ot H^0(S^{r-2}U^*)  \lra  V^*\ot (V^*\ot H^0(S^{r-2}U^*)) \lra  
V^*\ot H^0(S^{r-1}U^*).$$
This implies that a component of $H^0(S^{r-2}U^*)$ with length greater than
two on some factor will give components  of $V^*\ot H^0(S^{r-2}U^*)$ with the same property, 
and these cannot contribute to $D$.  \qed

\section{Equations of $\s_k(Seg(\BP A\times\BP B\times \BP C))$}

In this section we analyze the equations of  the secant varieties of 
$Seg(\BP A^*\times\BP B^*\times \BP C^*)$ for low dimensional vector spaces $A$, $B$, 
$C$.

For the dimensions of the secant varieties and filling $k$
that we use in this subsection, we refer the reader to \cite{CGG,Lickteig}.

\subsection{Case of $X=Seg(\pp m\times\pp n)$}\label{casemn}
\hspace*{3cm}\hfill \hfill\linebreak\indent 
Recall that here $\tdim \s_k(X)=k(m+n+2-k)-1$ until it fills and
$I(\s_k(X))$ is generated in degree $k+1$ by $\La {k+1}\BC^{m+1}
\ot\La {k+1}\BC^{n+1}$.

\subsection{Case of $X=Seg(\pp 1\times\pp 1\times\pp 1)$}\label{case111}
\hspace*{3cm}\hfill \hfill\linebreak\indent 
Here $\s_2(X)=\BP V$ and thus its ideal is zero.

\subsection{Case of $X=Seg(\pp 1\times\pp 1\times\pp c)$, $c=2,3$}
\label{case112}\label{case113}\hspace*{3cm}\hfill \hfill\linebreak
\indent  Here $I(\s_2(X))$ is generated in degree three  by flattenings 
and $\s_3(X)=\BP V$.  

\subsection{Case of $X=Seg(\pp 1\times\pp 2\times\pp 2)$ }
\label{case122} \hspace*{3cm}\hfill \hfill
 \linebreak\indent 
Again $I(\s_2(X))$ is generated in degree three by flattenings and $\s_3(X)=\BP V$.  
\medskip 
 
\subsection{Case of $X=Seg(\pp 2\times\pp 2\times \pp 2)$}\label{case222}
\hspace*{3cm}\hfill \hfill\linebreak\indent 
Again $I(\s_2(X))$ is generated in degree three by flattenings.  

\begin{proposition} Let $X=Seg(\pp 2\times\pp 2\times \pp 2)= 
Seg(\BP A^*\times \BP B^*\times\BP C^*)$.
\begin{itemize}
\item
The space of quartic equations of $\s_3(X)$ is
$$I_4(\s_3(X))
=S_{211}A\ot S_{211}B\ot S_{211}C, 
$$
and has dimension $27$. 

\item
The hypersurface $\s_4(X)$
is of degree nine and corresponds to the one-dimensional module 
$S_{333}A\ot S_{333}B\ot S_{333}C$.

\end{itemize}
\end{proposition}

A determinantal representation of these equations was given by Strassen, see \cite{GSS}.
We don't know if $I(\s_3(X))$ is generated in degree four.

This case is discussed in \cite{GSS} (without proofs).
To study $I_4(\s_3(X))$ we need only look at terms
$S_{\l_1}A\ot S_{\l_2} B\ot S_{\l_3}C$ with each  $\l_j$ of
length $3$ by case \ref{case122} since otherwise, by inheritance (proposition \ref{elimprop}), 
 we would have a nonzero element in $I_4(\s_3(\pp 1\times\pp 2\times\pp 2))$.
 Examining the decomposition
of $S^4(A\ot B\ot C)$ the only possible term is
$W=S_{211}A\ot S_{211}B\ot S_{211}C$, which occurs with multiplicity one.

To illustrate our methods, we give three proofs that $I_4(\s_3(Seg(\pp 2\times\pp 2\times \pp 2)))
=W$.

\smallskip

{\it First proof}:
We apply Proposition \ref{turbostuffin},
 that is, we   check   that $W$ is not contained in $V^2\ot S^2(V)$. 
This is easy: each term in $S^2V$ must have at least one symmetric power, say $S^2A$. If 
we tensor by the other $S^2A$ coming from $V^2$, we do not get the $S_{211}A$
term of $W$.  

\smallskip

{\it Second proof}: 
We make explicit  the embedding of $W^*$ in $S^4(A\ot B\ot C)$,
using the  first algorithm explained in \S\ref{schursect}. The representation $[211]$ of $\cS_4$ 
is the tensor product
 of the natural three dimensional 
representation $[31]$, given by the natural action of $\cS_4$ on the hyperplane
$x_1+x_2+x_3+x_4=0$ in $\CC^4$,
 with the sign representation. We choose the basis $e=(1,-1,0,0)$,  $f=(0,1,-1,0)$,  
$g=(0,0,1,-1)$. We keep this basis for $[211]$ but recall that the action of the symmetric group 
must be twisted by the sign. Averaging $e\ot e\ot f$ over the symmetric group, we obtain a non 
zero invariant $I$ in $[211]\ot [211]\ot [211]$, 
$$\begin{array}{rcl}
I & = & e\ot e\ot (-e+f+g)-f\ot f\ot (e-f+g)+g\ot g\ot (e+f-g) \\
 & & +(f-e)\ot (f-e)\ot (-e-f+g)+(g-e)\ot (g-e)\ot (e-f+g) \\
 & & +(g-f)\ot (g-f)\ot (e-f-g).
\end{array}$$
Now we need to evaluate the corresponding space of polynomials on $\s_3(X)$, which means 
that we may suppose, for example, that the third and fourth arguments are equal decomposed
tensors. Let $s_3$ denote the simple transposition $(34)$.
Since $s_3(e)=-e$, the contribution of all terms involving $e$ will vanish.
Moreover, $g=-s_3(f)$, and since the 
contributions of $f$ and $s_3(f)$ in our evaluation are obviously the same, $f$ and 
$-g$ have the same contribution. But if we let $e=0$ and $g=_f$ in the expression of 
the invariant $I$, we get zero, which means that our evaluation on $\s_3(X)$ does vanish. 

\smallskip 
{\it Third proof}: We use the second algorithm of section  \S\ref{schursect}. Taking the polynomial $P$
of example \ref{211211211ex} we see that if any two vectors are equal, each term of
$P$ vanishes.

\medskip
We now study   $\s_4(X)$. Since it is a hypersurface by \cite{CGG,Lickteig},
we need to look for instances of the
trivial representation in $S^d(A\ot B\ot C)$.
The first candidate appears when $d=6$, since $S^6(A\ot B\ot C)$ contains
$S_{222}A\ot S_{222}B\ot S_{222}C$.
The second candidate appears when $d=9$, since $S^9(A\ot B\ot C)$ contains
$S_{333}A\ot S_{333}B\ot S_{333}C$, again with multiplicity one. It is claimed in 
\cite{GSS} that the degree nine equation is the equation 
of $\s_4(X)$. We   verify this by applying Proposition \ref{turbostuffin}. 
We need  to see if the one dimensional representation $S_{333}A\ot S_{333}B\ot S_{333}C$ 
occurs in either $S^2V^3\ot V^2\ot V$ or $V^3\ot S^3V^2$. In the first  factor
of the first module, at least one of $A,B,C$,
say $A$, must occur as $S^2(S^3A)=S_{6}A\ot S_{42}A$. Since the partitions $(6)$ and $(42)$ 
are not contained in $(333)$, we cannot get $S_{333}A$ after tensoring by $V^2\ot V$.
For the second module, we note that $S_{333}A\ot S_{333}B\ot S_{333}C$ could only 
come from a factor $S_{33}A\ot S_{33}B\ot S_{33}C$ of $S^3V^2$. But 
$S^3(P\ot Q\ot R)$ does not contain 
$S_{111}P\ot S_{111}Q\ot S_{111}R$, so that up to symmetry,  either $S^3(S^2A)$ or 
$S_{21}(S^2A)$ must occur in each factor, and none of these contains $S_{33}A$.

\smallskip
Now that we have a nonzero  polynomial of degree nine that vanishes on the invariant 
hypersurface $\s_4(X)$, we conclude that it must be the equation of this hypersurface.
Indeed, suppose not. Then our polynomial would be the product of two polynomials, 
which automatically would be both invariant. In particular, their degrees would be 
multiples of three, so one of them would have degree three. But there is no invariant 
cubic polynomial.

The polynomial is described explicitly in example \ref{333333333ex} and one can verify that
it does indeed vanish on $\s_4(X)$, but some care must be taken in keeping track of
the signs. Similarly, one can explicitly write out the polynomial in
$S_{222}A\ot S_{222}B\ot S_{222}C$ and see that it does not vanish on $\s_4(X)$
(for example, this is easy to verify with Maple).

\smallskip

\subsection{Case of $X=Seg(\pp 1\times \pp 2\times\pp 3)$}
\label{case123}
\hspace*{3cm}\hfill \hfill

\begin{proposition} Let $X= Seg(\pp 1\times \pp 2\times \pp 3)= Seg(\BP A^*\times \BP B^*\times
\BP C^*)$. Then
$$
I_4(\s_3(X))
= S_{31}A\ot S_{211}B\ot S_{1111}C \op
S_{22}A\ot S_{22}B\ot S_{1111}C  
$$
is of dimension $6$ and   generates $I(\s_3(X))$.

Also,   $\s_4(X)=\BP V$.
\end{proposition}

\begin{proof}  By \cite{CGG,Lickteig},
$\tdim\s_3(X)=20$. On the other hand, $X\subseteq Seg(\pp 5\times\pp 3)$
and $\tdim\s_3(Seg(\pp 5\times\pp 3))=20$ (see case \ref{casemn} above).
Since both are irreducible varieties, they are equal and
$I_4(\s_3(X))=\La 4(A\ot B)\ot \La 4 C$.
 \end{proof}

\smallskip
 
\subsection{Case of $X=Seg(\pp 2\times \pp 2\times\pp 3)$}
\label{case223}
\hspace*{3cm}\hfill \hfill

\begin{proposition} Let $X= Seg(\pp 2\times \pp 2\times\pp 3)= Seg(\BP A^*\times \BP B^*\times
\BP C^*)$ . Then
\begin{itemize}

\item The space of quartics on $X$ is
\begin{align*}
I_4(\s_3(X))= &
\; S_{211}A\ot S_{211}B\ot S_{211}C \op 
S_{31}A\ot S_{211}B\ot S_{1111}C \\
&
\op S_{211}A\ot S_{31}B\ot S_{1111}C \op
S_{22}A\ot S_{22}B\ot S_{1111}C 
\end{align*}
 and is of dimension $135+2\times 45+36= 261$.

\item  $I_5(\s_4(X))=0$.
  
\item  $I_6(\s_4(X))= S_{222}A\ot S_{222}B\ot S_{3111}C \op S_{321}A\ot S_{321}B\ot S_{3111}C$ 
\newline
 The
  second term occurs with multiplicity one in $I_6$ while occuring with multiplicity
  four in $S^6V$. Note that $\tdim I_6(\s_4(X))=260$
    and it does not generate
  $I(\s_4(X))$ because the inherited $S_{333}A\ot S_{333}B\ot S_{333}C$ term in $I_9(\s_4(X))$
  cannot come from these terms.
  
\item   $\s_5(X)=\BP V$ 
\end{itemize}
 
\end{proposition}

We do not know if $I(\s_3(X))$ is generated in degree four.

\smallskip

To determine $I_4(\s_3(X))$,
in addition to $\La 4(A\ot B)\ot \La 4 C$,
  we inherit
$S_{211}A\ot S_{211}B\ot S_{211}C$ from   case \ref{case222}.
No other terms are possible by the same argument as in
case \ref{case222}.  

To see $I_5(\s_4(X))$ is empty we explicitly wrote down highest weight vectors in
all the possible modules and tested them at random points of $\s_4(X)$ with Maple.
We used the same method to find the modules in and not in $I_6(\s_5(X))$, but when we
found a polynomial that vanished, we checked the result symbolically.

In example \ref{3213213111} we gave an explicit basis of the highest weight
vectors for $S_{321}A\ot S_{321}B\ot S_{3111}C$. The linear combination
that vanishes  on $\s_4(X)$ is the polynomial obtained by symmetrizing 
$6F_1 -F_2 -4F_3+5F_4$.

The last assertion is not   in \cite{CGG,Lickteig} so we present a proof:

\begin{proof} We use Terracini's lemma.
Let $e_1\hd e_3$, $f_1\hd f_3$, $g_1\hd g_4$ respectively denote
bases of $\BC^3,\BC^3,\BC^4$. For our five points on $X$, take
$e_1\ot f_1\ot g_1$, $e_2\ot f_2\ot g_2$, $e_3\ot f_3\ot g_3$,
$(e_1+e_2+e_3)\ot (f_1+f_2+f_3)\ot g_4$,
$(e_1+e_2+e_3)\ot (f_1+\a f_2+\b f_3)\ot (g_1+g_2+g_3+ g_4)$
where $\a,\b$ are relatively prime and $|\a-\b|\neq 1$.
An easy calculation shows that if we use the monomial basis except for
using 
$e_1\ot f_2\ot g_3 +e_2\ot f_1\ot g_3$,
$e_2\ot f_3\ot g_1 +e_3\ot f_2\ot g_1$,
$e_1\ot f_3\ot g_2 +e_3\ot f_1\ot g_2$,
$e_1\ot f_3\ot g_4 +e_3\ot f_1\ot g_4$,
$e_1\ot f_2\ot g_3 -e_2\ot f_1\ot g_3$,
$e_2\ot f_3\ot g_1 -e_3\ot f_2\ot g_1$,
$e_1\ot f_3\ot g_2 -e_3\ot f_1\ot g_2$,
$e_1\ot f_3\ot g_4 -e_3\ot f_1\ot g_4$, 
instead of the monomials that appear in them, then the span of
the tangent spaces to the first four points is all but the last four
terms, and adding tangent space to the fifth point enables us
to dispense with those.
(Recall that $\hat T_{[e\ot f\ot g]}
Seg(\BP E\times \BP F\times \BP G)=
E\ot f\ot g + e\ot F\ot g + e\ot f\ot G$.)
\end{proof}

\smallskip

\subsection{Case of $X=Seg(\pp 2\times \pp 3\times\pp 3)$}\label{case233}
\hspace*{3cm}\hfill \hfill

\begin{proposition}  Let $X= Seg(\pp 2\times \pp 3\times\pp 3)= Seg(\BP A^*\times \BP B^*\times
\BP C^*)$ . Then
\begin{itemize} 

\item The space of quartics on $\s_3(X)$ is
\begin{align*}
I_4(\s_3(X))=&
\quad S_{211}A\ot S_{31}B\ot S_{1111}C  \op S_{211}A\ot S_{1111}B\ot S_{31}C \\
& \quad\op S_{22}A\ot S_{22}B\ot S_{1111}C\op S_{22}A\ot S_{1111}B\ot S_{22}C  \\
 & \quad\op S_{31}A\ot S_{1111}B\ot S_{211}C   \op S_{31}A\ot S_{211}B\ot S_{1111}C \\   
& \quad\op S_{4}A\ot S_{1111}B\ot S_{1111}C \op S_{211}A\ot S_{211}B\ot S_{211}C
\end{align*}
and has  dimension $2\times 135+2\times 120+2\times 225+15+675=1650$. 

\item The space of quintic equations of $\s_4(X)$ is 
$$I_5(\s_4(X))=S_{311}A\ot S_{2111}B\ot S_{2111}C, $$
and has dimension $96$. $I(\s_4(X))$ is not generated in degree five.

\item  $I_6(\s_5(X))
=0
$. 

\item $I_7(\s_5(X))
=0
$. 
 
\end{itemize}
\end{proposition}

Note that   $\s_6(X)$ fills, see \cite{CGG,Lickteig}.  

We do not know if the ideal of $\s_3(X)$ is generated in degree four.

\smallskip

\begin{proof}

$I_4(\s_3(X))$ follows from flattenings and inheritance.

Continuing to $\s_4(X)$,
since $I_5(\s_4(Seg(\pp 2\times \pp 2\times\pp 3))=0$, 
the only possible term in $I_5(\s_4(X))$ is   $W=S_{311}A\ot S_{2111}B\ot S_{2111}C$,
which occurs in $S^5(A\ot B\ot C) $ with multiplicity one, because this is the unique 
component with a partition of length three in the $A$ factor and length four in the $B,C$ 
factors.  But this factor does not occur inside $S^3V\ot V^2$.
Thus Proposition \ref{turbostuffin} applies.

Turning to $\s_5(X)$, since $\s_5$ fills for both
$\pp 2\times \pp 2\times \pp 3$ and $\pp 1\times\pp 3\times \pp 3$
we only need look at elements of $S^6(A\ot B\ot C)$ that are of length
three in the first factor and four in the second and third factors.
Examining the decomposition, the candidate modules
(up to permutation in the last two factors)  are:
$S_{411}A\ot S_{3111}B\ot S_{3111}C$, $S_{411}A\ot S_{2211}B\ot S_{2211}C$,
 $S_{222}A\ot S_{3111}B\ot S_{3111}C$, $S_{222}A\ot S_{2211}B\ot S_{2211}C$
 with multiplicity one and
$S_{411}A\ot S_{3111}B\ot S_{2211}C$,
$S_{321}A\ot S_{3111}B\ot \linebreak S_{3111}C$,
$S_{321}A\ot S_{3111}B\ot S_{2211}C$, $S_{321}A\ot S_{2211}B\ot S_{2211}C$ 
with multiplicity two.

On the other hand, consider $S^4V\ot V^2$. In order to have
two modules with partition of length four, we need the partitions
in $S^4V$ to have length at least three. All are accounted for,
so Proposition \ref{turbostuffin} is not useful here.
Thus we do direct calculations with Maple, which is what we use for
  $I_7(\s_5(X))$  as well, the latter being quite involved as modules
 appear with multiplicity up to nine.
 \end{proof}

 \smallskip

\subsection{Case of $X=Seg(\pp 3\times \pp 3\times\pp 3)$}\label{case333}
\hspace*{3cm}\hfill \hfill

\begin{proposition} Let $X= Seg(\pp 3\times \pp 3\times\pp 3)= Seg(\BP A^*\times \BP B^*\times
\BP C^*)$ . Then
\begin{itemize}

\item $I_4(\s_3(X)) = \La 4(A\ot B)\ot \La 4 C 
$ plus permutations and
$S_{211}A\ot S_{211}B\ot S_{211}C$.

\item The space of quintic equations of $\s_4(X)$ is 
$$I_5(\s_4(X))=S_{311}A\ot S_{2111}B\ot S_{2111}C\op S_{2111}A\ot S_{311}B\ot S_{2111}C\op
S_{2111}A\ot S_{2111}B\ot S_{311}C$$
and has dimension $3\times 36\times 4\times 4=1728$.  

 We also know that $S_{333}A\ot S_{333}B\ot S_{333}C$
is in $I_9(\s_4(X))$ by inheritance. Since it only involves partitions of length three,
it cannot be generated by $I_5(\s_4(X))$, whose components all involve partitions 
of length four.   Thus $I(\s_4(X))$ is not generated in degree $5$.  
  
 \item 
$ 
I_6(\s_5(X)) =0$

\item $ 
I_7(\s_5(X)) =0$

\item $I_8(\s_5(X))
\supseteq S_{5111}A\ot S_{2222}B\ot S_{2222}C \op S_{3311}A\ot S_{2222}B\ot S_{2222}C,
$  again, up to permutations. Both modules occur with multiplicity one in $S^8(A\ot B\ot C)$.

\item $I_d(\s_6(X)) =0$ for $d\leq 8$.

\end{itemize}
\end{proposition}
 
The factors $S_{311}A\ot S_{2111}B\ot S_{2111}C$ plus permutations
in $I_5(\s_4(X))$
are inherited from case \ref{case233}.
Since $S_{2111}A\ot S_{2111}B\ot S_{2111}C$ is not in $S^5V$, all of
$I_5(\s_4(X))$ must be inherited from case \ref{case233}. 

The remaining modules were eliminated by extensive Maple calculations.
These calculations also showed us the candidate members of $I_8$ but
only with extremely high probability, so we now present direct proofs
that they are in the ideal.

 The following monomial gives a 
 highest weight vector for $S_{5111}A\ot S_{2222}B\ot S_{2222}C$  when summed over the symmetric group:
$$F =  \a_1\a_2\a_5\a_6(\a_3\we \a_4\we \a_7\we \a_8)(\b_1\we \b_2\we \b_3\we \b_8)(\b_4\we \b_5\we 
\b_6\we \b_7)
(\g_1\we \g_2\we \g_3\we \g_4)(\g_5\we \g_6\we \g_7\we \g_8). 
$$

A general element of $\s_5(\pp 3\times\pp 3\times \pp 3)$ is of the form $a_1b_1c_1+\cdots +a_5b_5c_5$, 
and when we compute a homogeneous polynomial $P$ on such a sum, we get, after expansion, terms with different 
homogeneities on $a_1b_1c_1, \ldots ,a_5b_5c_5$. These homogeneous components must all vanish identically
if we want $P$ to vanish on $\s_5(\pp 3\times\pp 3\times \pp 3)$. Note that the type of 
homogeneity, up to permutation of $1...5$, is given by a partition $\pi$ with $5$ parts, 
the sum of the parts being equal to the degree of the polynomial $P$. 

\smallskip
We associate a graph  $\g (F )$   to our 
tensor $F$. The vertices are identified with the integers $1...8$, and two 
vertices $i$ and $j$ are joined by an edge iff they are not wedged together in the 
expression of $F$. We get:\bigskip

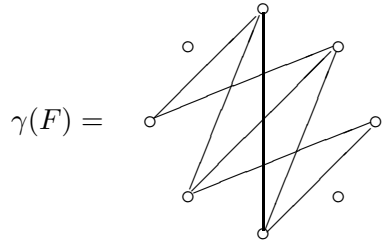
\begin{center}
\setlength{\unitlength}{5mm}
\begin{picture}(20,6)(-7,0)
\put(-3.5,3){$\gamma (F)=$}
\multiput(0,3)(6,0){2}{$\circ$}
\multiput(3,0)(0,6){2}{$\circ$}
\multiput(1,1)(4,0){2}{$\circ$}
\multiput(1,5)(4,0){2}{$\circ$}
\put(0.3,3.3){\line(1,1){2.7}}
\put(1.25,1.4){\line(2,5){1.85}}
\put(3.2,0.3){\line(0,1){5.8}}
\put(0.3,3.3){\line(5,2){4.75}}
\put(1.3,1.35){\line(1,1){3.7}}
\put(3.25,0.3){\line(2,5){1.9}}
\put(1.35,1.3){\line(5,2){4.7}}
\put(3.3,0.25){\line(1,1){2.75}}
\end{picture}
\end{center}

Observe that $\gamma (F )$ contains {\it no triangle}, 
so that each time we choose a triple of indices among $1...8$, two of 
them are wedged together somewhere in the expression of $F$. 

Thus when we evaluate $P$ on the monomials in the expansion of
$(a_1b_1c_1+\cdots +a_5b_5c_5)^5$, all terms with a power of three or greater
evaluate to zero.
 
There remains to consider the case  where the degrees are 
$(2,2,2,1,1)$ (the case of $(2,2,2,2)$ will follow).
Denote the indices occuring with a power $2$ by $s,t,u$ and
those to the first power by $i,j$.
Note that $s,t,u$ must appear twice in the contributions 
of $A$, $B$, $C$, but of course not in a same wedge product. So we'll only get terms
of type 
$$\a_{s}\a_{t}\a_{\g}\a_i(\a_{s}\we \a_{t}\we \a_{u}\we \a_j)
(\b_{s}\we \b_{t}\we \b_{u}\we \b_i)
(\b_{s}\we \b_{t}\we \b_{u}\we \b_j)
(\g_{s}\we \g_{t}\we \g_{u}\we \g_i)(\g_{s}\we \g_{t}\we \g_{u}\we \g_j).$$
But this is skew-symmetric, e.g., in $s$ and $t$, so the total contribution 
of these kinds of terms is zero. \medskip


For $S_{3311}A\ot S_{2222}B\ot S_{2222}C=S_{22}A\ot \det A\ot (\det B)^2\ot (\det C)^2$  
the analysis is similar to the previous case.  Here we may take
$$   F    =  (\a_1\we \a_3)(\a_5\we \a_7)(\a_2\we \a_4\we \a_6\we \a_8)(\b_1\we \b_2\we \b_5\we \b_6)
(\b_3\we \b_4\we \b_7\we \b_8)(\g_1\we \g_2\we \g_3\we \g_4)(\g_5\we \g_6\we \g_7\we \g_8).
$$

The associated graph is as follows. Again, it contains no triangle:
\bigskip

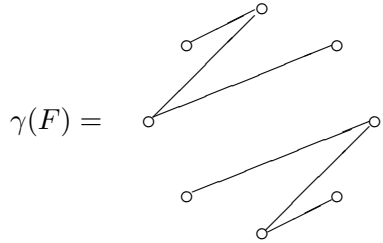
\begin{center}
\setlength{\unitlength}{5mm}
\begin{picture}(20,6)(-1,0)
\put(3,3){$\gamma (F)=$}
\multiput(6.5,3)(6,0){2}{$\circ$}
\multiput(9.5,0)(0,6){2}{$\circ$}
\multiput(7.5,1)(4,0){2}{$\circ$}
\multiput(7.5,5)(4,0){2}{$\circ$}
\put(6.8,3.3){\line(1,1){2.7}}
\put(7.8,5.3){\line(2,1){1.7}}
\put(6.8,3.3){\line(5,2){4.75}}
\put(7.85,1.3){\line(5,2){4.7}}
\put(9.8,0.3){\line(1,1){2.75}}
\put(9.85,0.25){\line(2,1){1.7}}
\end{picture}
\end{center}

Finally,  
consider the   terms we   get with exponents $ (2,2,2,1,1)$. 
They must be of type 
$$(\a_{s}\we \a_{t})(\a_{u}\we \a_i)(\a_{\a}\we \a_{t}\we \a_{u}\we \a_j)
(\b_{s}\we \b_{t}\we \b_{u}\we \b_i)
(\b_{s}\we \b_{t}\we \b_{u}\we \b_j)
(\g_{s}\we \g_{t}\we \g_{u}\we \g_i)(\g_{s}\we \g_{t}\we \g_{u}\we \g_j).$$
This is no longer skew-symmetric in $s,t$. But if we symmetrize with 
respect to $s,t,u$, we get (twice) the product of a fixed product of 
determinants, with $(\a_{s}\we \a_{t})(\a_{u}\we \a_i)+(\a_{t}\we \a_{u})(\a_{s}\we \a_i)+
(\a_{u}\we \a_{s})(\a_{t}\we \a_i)$. Since the vanishing of such an expression is precisely 
the condition that defines $S_{22}A$ inside $S^2(\wedge^2A)$, our proof is complete. 

In degree nine we verified that all cases of low multiplicity do not
arise in $I(\s_6(X))$ and we are currently working on the cases of
higher multiplicity. However, inspired by the $\pp 2\times \pp 2\times \pp 2$ case,
there are two natural candidates that we checked using $C++$ code
written by P. Barbe:

\begin{proposition}
The module  

$$S_{3333}A\ot S_{3333}B\ot S_{3333}C   \subset S^{12}(A\ot B\ot C),$$
  which occurs with multiplicity one,
 is not in $I(\s_6(\pp 3\times\pp 3\times \pp 3))$.
\end{proposition}

The polynomial in degree $12$ may be obtained by symmetrizing
\begin{align*}
F(a_1b_1c_1\hd a_{12}b_{12}c_{12}) = & 
det(a_1,a_2,a_3,a_4)det(a_5,a_6,a_7,a_8)det(a_9,a_{10},a_{11},a_{12})\\
 & 
 det(b_1,b_2,b_5,b_6)det(b_3,b_7,b_9,b_{10})det(b_4,b_8,b_{11},b_{12})\\
& 
det(c_1,c_7,c_9,c_{12})det(c_3,c_5,c_8,c_{10})det(c_2,c_4,c_6,c_{11}).
\end{align*}
 



\begin{thebibliography}{aa}

\bibitem{AH} Alexander J., Hirschowitz A., {\sl
 Polynomial interpolation in several variables},  J.
   Algebraic Geom. {\bf 4}  (1995), no. 2, 201--222.
   
   \bibitem{bou} Bourbaki N., Groupes et alg\`ebres de Lie, Hermann,
Paris 1968. 

\bibitem{CGG} Catalisano M. V., Geramita A.V., Gimigliano
 A.,  {\sl Ranks of tensors, secant varieties
   of Segre varieties and fat points},
    Linear Algebra Appl. {\bf 355} (2002), 263--285.
    
\bibitem{CS} Gonzalo Comas, Malena Seiguer, {\sl On the rank of a binary form},
    arXiv:math.AG/0112311.
    
 \bibitem{FH} Fulton, W. and Harris J. {\sl Representation Theory, a first course},
 Springer Verlag, GTM 129 (1991).
  
    \bibitem{GSS} Luis David Garcia, Michael Stillman, Bernd Sturmfels, 
  {\sl Algebraic Geometry of Bayesian Networks}, arXiv:math.AG/0301255.
  
  

\bibitem{GH} Griffiths P.A., Harris J., {\it Algebraic Geometry and
Local Differential Geometry }, Ann. scient. Ec. Norm. Sup.  
{\bf 12}, 355-432 (1979).
  
  \bibitem{Lsec} Landsberg J.M., {\it On degenerate secant and tangential
varieties and local differential  geometry},  Duke Math. J. {\bf 85}  (1996),
605--634.
  
  \bibitem{LM0} Landsberg J.M., Manivel L., {\sl On the projective
geometry of homogeneous varieties}, Comm. Math. Helv. {\bf 78} (2003), 65--100.

 \bibitem{LMuniv} Landsberg J.M., Manivel L., {\sl A universal dimension formula
for complex simple Lie algebras}, preprint.

\bibitem{Lickteig}
Lickteig, Thomas {\it Typical tensorial rank},   Linear Algebra Appl.  {\bf 69}  (1985), 95--120.


 
 \bibitem{Strassenb} Strassen V., {\sl Rank and optimal computation of generic tensors}, 
 Linear algebra Appl.   {\bf 52/53} (1983), 645-685.
 
  \bibitem{Strassen} Strassen V., {\sl The asymptotic spectrum of tensors}, 
 Crelles J. Reine. Angew. Math. {\bf 384} (1988), 102-152.


\bibitem{weyman} Weyman J., Cohomology of vector bundles and syzygies, Cambridge
Tracts in Mathematics 149, Cambridge University Press 2003. 

  
\bibitem{Winograd} Winograd, S.
{\it On multiplication of $2\times 2$ matrices},
Linear Algebra and Appl. {\bf 4} (1971), 381--388.


\bibitem{Zak} Zak F.: Tangents and secants of algebraic varieties, 
Translations of Math. Monographs {\bf 127}, AMS 1993. 

\end{thebibliography}
\end{document}